\documentclass[a4paper,10pt]{article}

\usepackage{amsfonts,amssymb,mathrsfs,amscd}
\def \R {{\mathbb {R}}}

\usepackage[T1]{fontenc} 
\usepackage[cp1251]{inputenc}  
\usepackage[russian]{babel}

\begin{document}
\title{\bf Singular geometric averages \\ for ergodic multiflows}
\author{\bf I.V. Bychkov, V.V. Ryzhikov }
\date{}

\maketitle
\Large

We consider ergodic actions of groups $\R^n$ (multiflows) on a probability space. The general theorem on universal averaging is applied to averaging along manifolds. A class of ergodic theorems arises in which singular averaging is performed over a sequence of expanding smooth manifolds. For an ergodic multiflow, convergence of time averages holds, for example, over sequences of expanding spheres or smooth curves in general position.

For ergodic flows on the probability space $(X,\mu)$, Kozlov and Treshchev in \cite{KT} considered non-uniform averages of the form
$$P_tf(x)=\int_R f(T_{rt}x) h(r) dr, \ \int_R h\,dr=1,\ h\geq 0.$$
In particular, in \cite{KT}, the convergence of $P_tf$ to a constant as $t\to +\infty$ was proved for functions $f\L_1$, which strengthens von Neumann's theorem, in which $h=\chi_{[0,1]}$. We consider ergodic flows $T_t$, with multidimensional time $t\in \R^d$. The ergodicity property for a multiflow $T_t$ that preserves the measure $\mu$ means that every $T_t$-invariant $\mu$-measurable function is constant.

The method for proving the aforementioned Kozlov-Treshchev result leads, without fundamental changes, to the following generalization.

\vspace{2mm}
\bf Theorem 1. \it  Let $\nu$ be a normalized measure on $\R^d$ that is absolutely continuous with respect to Lebesgue measure. Then, for an ergodic multiflow $T_t$ that preserves the probability measure $\mu$, and
a function $f\in L_1(\mu)$, the convergence holds:
$$\left\|P_tf-\int f \, d\mu\right\|_1\to 0, \ t\to +\infty,$$ 
where
$$P_tf(x)=\int_{R^d} f(T_{tr}x) d\nu(r),\ t\in \R.$$
\rm

\vspace{2mm}
\bf Universal averaging and convolution of weights. \rm Following \cite{R}, a normalized Borel measure $\nu$ on $R^d$ will be called \it universal if for every ergodic multiflow $T_t$ and every function $ f\in L_1(\mu)$ the statement of Theorem 1 holds.

\vspace{2mm}
\bf Theorem 2. \it If some convolution power $\nu^{\ast n}$ of a normalized Borel measure $\nu$ on $\R^d$ is universal for an ergodic flow, then the measure $\nu$ is also universal (corollary of Theorem 3). \rm

\vspace{2mm}
A sequence of normalized Borel measures $\nu_j$ on a locally compact commutative group $G$ is called a universal sequence if
for every ergodic action $\{T_g\,:\, g\in G\}$ and every function $f\in L_1(\mu)$,
$$\left\|P_jf -\int f \, d\mu\right\|_1\to 0, \ \ j\to +\infty,$$
where $P_j$ is the operator
defined by $P_jf(x)=\int_{G} f(T_{g}x) d\nu_j(g).$

\vspace{2mm}
\bf Theorem 3. \it Let $\nu_j$ be a sequence of normalized Borel measures on a locally compact commutative group $G$. If for some $n$ the convolution powers $\nu_j^{\ast n}$ form a universal sequence for ergodic actions of $G$, then the sequence $\nu_j$ is also universal. \rm

\vspace{2mm}
The proof of this theorem is similar to the case $G=\R$ \cite{R}. The scheme is as follows. If $P_j$ is the averaging operator with respect to the measure $\nu_j$, then $P^n_j$ is the operator corresponding to the weight $\nu^{\ast n}_j=\nu_j\ast \dots\ast \nu$ ($n$ factors). For normal operators $P_j$ and a function $f\in L_2$, the convergence $$\left\|P^n_jf -\int f d\mu\right\|_2 \to 0$$ is equivalent to the convergence $$\left\|P^n_jf -\int f d\mu\right\|_2 \to 0.$$ This readily implies the convergence of $P_jf$ to a constant in $L_1$. It remains to use the density of $L_2$ in $L_1$.

\newpage
\bf Geometric Singular Averages. \rm
From Theorems 2 and 3, we can obtain a variety of statements about the convergence of means over expanding manifolds in $\R^d$, $d>1$, of dimension less than $d$. We formulate two such corollaries.

\vspace{2mm}
\bf Theorem 4. \it  Let the normalized measure $\nu$ be absolutely continuous with respect to Lebesgue measure on a smooth curve $\Gamma\subset \R^d$ in general position, i.e., for almost all sets of points $t_1,\dots, t_d\in\Gamma$, the tangent vectors at these points form a linearly independent system. Then the convolution degree $\nu^{\ast d}$ is absolutely continuous with respect to the Lebesgue measure in $\R^d$, hence the measure $\nu$ is universal. \rm

\vspace{2mm}
Let's give a concrete example. In the space $\R^d$, consider
the curve $$\Gamma= \{(r, r^2,\dots, r^d)\,:\, 0< r < 1\}.$$
The vector $(1,2r, \dots, (d-1)r^{d-1}$ is collinear with the tangent vector of our curve at the point $(r, r^2,\dots, r^d)$. The set of $d$ distinct
vectors belonging to the curve is linearly independent, which is a consequence of the non-singularity of the Vandermonde matrix for sets of distinct non-zero numbers $r_1, r_2,\dots, r_{d}$.

Wiener's ergodic theorem \cite{W} considered averaging over balls. Spheres can be substituted for balls without losing the universality of the average.

\vspace{2mm}
\bf Theorem 5. \it Let $S_j$ be a sequence of spheres of dimension $d$ in the space $\R^{d+1}$.
If the radii of these spheres tend to infinity, then uniform averages over these spheres for the $(d+1)$-parametric ergodic flow converge to constant functions. \rm

\vspace{2mm}
Consider the case $d=2$. Let $\nu$ denote the normalized Lebesgue measure on the two-dimensional sphere $S\subset \R^3$, and $\sigma_\gamma$ the Lebesgue measure on its semimeridian passing through the equatorial point $\gamma$. Then the convolution $\sigma_\gamma\ast \nu$ in $\R^3$ is absolutely continuous with respect to the Lebesgue measure $m$ in $\R^3$. This follows from the fact that for almost all pairs of vectors of the form (meridian point, spherical point), the Jacobian of the map of this pair of vectors to their sum is nonzero. Now we note that the convolution $\nu\ast \nu$ is absolutely continuous with respect to the integral over the variable $\gamma$ of the family of measures $\sigma_\gamma\ast \nu$. The absolute continuity of these measures with respect to $m$ implies that $\nu\ast \nu$ is also absolutely continuous with respect to $m$. The proof of the general case uses a similar technique.

\vspace{3mm}
Wiener proved in \cite{W} the almost-everywhere convergence of spherical averages.
An interesting question arises about the almost-everywhere convergence
of spherical and more general averages over manifolds of dimension
less than $d$.  

\vspace{3mm}
Another question:
can  be universal a normalized measure on $\R^d$ such that all its convolution powers are  singular?

 \vspace{12mm}
\ \ \ \ \ \ \ \ \ \ \ \ \ \ \ \ \ \ \ \ \ \ \ \bf References  \rm

 \vspace{3mm}
[1]\ V.V. Kozlov, D.V. Treschev, On new forms of the ergodic theorem, J. Dynam. Control Systems, 9:3 (2003), 449-453

 \vspace{3mm}
[2] \ V.V. Ryzhikov, Universal Averaging for Ergodic Flows, Mat. Notes, 119:5 (2026), 759-765

 \vspace{3mm}
[3] \	N. Wiener, The ergodic theorem, Duke Math. J., 5:1 (1939), 1-18

\newpage

{\bf \ \ \ \ \ \ \ \ \  \ Сингулярные геометрические усреднения   

\ \ \ \ \ \ \ \ \ \ \ \ для эргодических мультипотоков}

\vspace{5mm}
\bf \ \ \ \ \ \ \ \ \ \ \ \ \ \ \ \ И.В. Бычков, В.В. Рыжиков \rm

\vspace{10mm}
Рассматриваются эргодические действия групп $\R^n$ (мультипотоки) на вероятностном пространстве.   К усреднениям вдоль многообразий  применяется общая теорема об универсальных  усреднениях.  Возникает   класс  эргодических теорем, в которых сингулярные усреднения осуществляются по последовательности  расширяющихся гладких многообразий.   Для эргодического мультипотока имеет место сходимость  временных средних, например, по последовательностям расширяющихся сфер, или гладких кривых общего положения.

\vspace{5mm}
Для эргодических потоков на вероятностном пространстве $(X,\mu)$ Козловым и Трещевым в \cite{KT}  рассматривались   временные неравномерные  усреднения вида  
$P_tf(x)=\int_R  f(T_{rt}x) h(r) dr,$ \ $\int_R h\,dr=1,$ \ $h\geq 0.$
В частности,  в \cite{KT}  была    для функций $f\L_1$ была доказана   сходимость $P_tf$  к константе    при $t\to +\infty$, что    усиливает теорему фон Неймана, в которой  $h=\chi_{[0,1]}$.  Рассматрим  эогодические потоки  $T_t$, с многомерным временем $t\in \R^d$.  Свойство эргодичности для мультипотока  $T_t$, сохраняющего меру $\mu$,   означает, что всякая $T_t$-инвариантная  $\mu$-измеримая функция является константой.  
Метод доказательства упомянутого результата Козлова-Трещева   без принципиальных изменений   приводит к следующему обобщению. 

 \vspace{2mm}
\bf Теорема 1. \it Пусть $\nu$ -- нормированная мера на $\R^d$, абсолютно непрерывная   относительно меры Лебега. Тогда для эргодического мультипотока $T_t$, сохраняющего вероятностую меру $\mu$, и
 функции   $f\in L_1(\mu)$ имеет место  сходимость  
$$\left\|P_tf-\int f \, d\mu\right\|_1\to 0,  \ t\to +\infty,$$ где
$$P_tf(x)=\int_{R^d}  f(T_{tr}x) d\nu(r),\ t\in \R.$$
\rm

\vspace{2mm}
 \bf Унивесальные усреднения и свертки весов. \rm Следуя  \cite{R},  нормированную борелевскую  меру $\nu$ на $R^d$ будем называть \it универсальной,  \rm  если   для всякого эргодического мультипотока $T_t$ и всякой функции   $ f\in L_1(\mu)$  выполнено утверждение теоремы 1. 

\vspace{2mm}
 \bf Теорема 2.  \it   Если  некоторая  сверточная степень $\nu^{\ast n}$   нормированной борелевской меры $\nu$ на $\R^d$  универсальна для эргодического потока, то мера $\nu$  также универсальна  (следствие теоремы 3).  \rm

\vspace{2mm}
Последовательность нормированных борелевских мер $\nu_j$  на локально компактной коммутативной группе $G$  называется \it универсальной  последовательностью, \rm  если
для всякого эргодического действия  $\{T_g\,:\, g\in G\}$ и  всякой функции $f\in L_1(\mu)$  выполнено 
$$\left\|P_jf -\int f \, d\mu\right\|_1\to 0, \ \ j\to +\infty,$$
где  $P_j$ -- оператор, 
определенный равенством $$P_jf(x)=\int_{G}  f(T_{g}x)  d\nu_j(g).$$

\vspace{2mm}
\bf Теорема 3.  \it  Пусть $\nu_j$ --  последовательность нормированных борелевских мер на локально компактной коммутативной группе $G$. Если  для некоторого $n$ сверточные степени $\nu_j^{\ast n}$ образуют  универсальную  последовательность для эргодических действий группы $G$, то последовательность $\nu_j$ также является  универсальной. \rm 

\vspace{2mm}
Доказательство этой теоремы проводится аналогично  случаю $G=\R$ \cite{R}. Схема следующая. Если $P_j$  -- оператор усреднения по мере $\nu_j$,  то $P^n_j$ --   оператор, отвечающий весу    $\nu^{\ast  n}_j=\nu_j\ast \dots\ast \nu$ ($n$ сомножителей).   Для нормальных операторов $P_j$ и  функции $f\in L_2$ сходимость $\left\|P^n_jf -\int f d\mu\right\|_2 \to 0$ эквивалентна  сходимости $\left\|P^n_jf -\int f d\mu\right\|_2 \to 0$.   Из этого несложным образом вытекает  сходимость $P_jf$ к константе  в  $L_1$. Остается воспользоваться плотностью   $L_2$ в $L_1$.

\newpage
\bf Геометрические сингулярные усреднения. \rm
  Из теорем 2, 3  можно получить  многообразие  утверждений о сходимости средних по расширяющимся многообразиям в $\R^d$, $d>1$, размерности меньшей, чем $d$. Сформулируют два таких следствия.

 \vspace{2mm}
\bf Теорема 4.  \it Пусть нормированная мера $\nu$   абсолютна непрерывна   относительно  меры Лебега на гладкой кривой $\Gamma\subset \R^d$ общего положения, т.е.   для  почти всех наборов точек $t_1,\dots, t_d\in\Gamma$ касательные векторы в этих точках  обрзуют линейно независимую систему. Тогда   сверточная степень $\nu^{\ast d}$ абсолютно непрерывна относительно  меры Лебега в $\R^d$, следовательно, мера $\nu$ универсальна. \rm

 \vspace{2mm}
Приведем конкретный  пример.     В пространстве $\R^d$ рассмотрим   
кривую   $$\Gamma= \{(r, r^2,\dots, r^d)\,:\, 0< r < 1\}.$$
Вектор  $(1,2r, \dots, (d-1)r^{d-1}$   коллинеарен  касательному вектору  нашей кривой в точке $(r, r^2,\dots, r^d)$.  Набор из $d$  различных 
векторов, принадлежащих кривой,  линейно независим, что является 
следствием невырождености матрицы Вандермонда  для наборов различных ненулевых чисел $r_1, r_2,\dots, r_{d}$.  

\vspace{2mm}
В эргодической теореме  Винера \cite{W} рассматривались  усреднения по шарам. Заменим шары   на сферы.  

\vspace{2mm}
\bf Теорема 5. \it Пусть $S_j$  последовательность сфер размерности $d$  в пространстве $\R^{d+1}$.
Если радиусы этих сфер стремятся к бесконечности, то равномерные  усреднения по этим сферам   для $(d+1)$-параметрического эргодического потока сходятся к постоянным функциям.  \rm

\vspace{2mm}
Рассмотрим случай $d=2$.   Пусть $\nu$ обозначает  нормированную меру Лебега на  двумерной сфере $S\subset \R^3$, а $\sigma_\gamma$ -- меру Лебега на ее  полумеридиане, проходящим через точку экватора  $\gamma$.  Тогда свертка $\sigma_\gamma\ast \nu$ в $\R^3$ абсолютно непрерывна относительно меры Лебега $m$  в $\R^3$. Это следует из того, что  для почти всех пар векторов вида (точка меридиана, точка сферы) будет ненулевым якобиан отображения этой пары векторов в их сумму.   Теперь замечаем, что  свертка  $\nu\ast \nu$  абсолютно  непрерывна относительно интеграла по    переменной $\gamma$ семейства   мер $\sigma_\gamma\ast \nu$.  Из абсолютной непрерывности этих мер относительно  $m$ вытекает, что  $\nu\ast \nu$  также  абсолютно непрерывна относительно  $m$. Доказательство общего  случая  использует  похожий прием.

\vspace{3mm}
Винер доказал в \cite{W}  сходимость почти всюду для шаровых усреднений.
Возникает интересный вопрос о сходимости почти всюду 
 сферических  и более общих усреднений по   многообразиям  размерности,
меньшей $d$.  

\vspace{3mm}
Другой вопрос:
может ли мера  на $\R^d$, у которой все сверточные степени  сингулярны, быть универсальной для эргодических действий группы $\R^d$?


\begin{thebibliography}{9}
\bibitem{KT} V.V. Kozlov, D.V. Treschev, On new forms of the ergodic theorem, J. Dynam. Control Systems, 9:3 (2003), 449-453


\bibitem{R}V.V. Ryzhikov, Universal Averaging for Ergodic Flows, Mat. Notes, 119:5 (2026), 759-765

\bibitem{W}	N. Wiener, The ergodic theorem, Duke Math. J., 5:1 (1939), 1-18

\end{thebibliography}
\end{document}